\documentclass[reqno]{amsproc}
\usepackage{amsmath}
\usepackage{amsfonts}
\usepackage{mathrsfs}
\usepackage{amssymb}
\usepackage{circuitikz}
\usepackage{color}
\usepackage{tikz}
\usepackage{needspace}
\usepackage[hidelinks]{hyperref}

\makeatletter
\@namedef{subjclassname@2010}{%
\textup{2010} Mathematics Subject
Classification} \makeatother

\DeclareMathOperator{\lin}{lin}

\DeclareMathOperator{\dist}{dist}

\numberwithin{equation}{section}

\newtheorem{theorem}{Theorem}[section]

\newtheorem{lemma}[theorem]{Lemma}

\theoremstyle{remark}
\newtheorem{df}[theorem]{Definition}

\newcommand*{\ascr}{\mathscr{A}}

\newcommand*{\hh}{\mathcal{H}}

\newcommand*{\natu}{\mathbb{N}}

\newcommand*{\comp}{\mathbb{C}}
\newcommand*{\borel}{\mathfrak{B}}

\newcommand*{\D}{\mathrm{d\hspace{.1ex}}}
\newcommand{\abs}[1]{\lvert #1\rvert}
\theoremstyle{definition}

\newcommand*{\Le}{\leqslant}
\newcommand*{\Ge}{\geqslant}

\newcommand*{\ogr}{\boldsymbol{B}}
\newcommand*{\rbb}{\mathbb{R}}

\newcommand{\norm}[1]{\lVert #1\rVert}

\title[Hyperrigidity and convexity in two variables]
{Hyperrigidity and convexity in two variables}

   \author[P. Pietrzycki]{Pawe{\l} Pietrzycki}

   \subjclass[2020]{Primary 46L07, 47A20, 47A63,
47B15; Secondary 44A60, 46G10}

   \keywords{Hyperrigidity, completely positive map,
dilations, (semi-) spectral measure, weak and strong operator topologies, Choquet
boundary}

   \address{Wydzia{\l} Matematyki i Informatyki, Uniwersytet
Jagiello\'{n}ski, ul. {\L}ojasiewicza 6, PL-30348
Krak\'{o}w, Poland}

   \email{pawel.pietrzycki@im.uj.edu.pl}

\begin{document}
\begin{abstract}
We extend Brown's one-variable rigidity theorem to strictly convex
functions of two variables under a regularity assumption. Let
$K\subset\rbb^2$ be nonempty, compact, and convex, and let $f$ be
real-valued, strictly convex on $K$, and of class $C^2$ on a neighbourhood
of $K$. We prove that the generating set $\{x_1,x_2,f\}$ is hyperrigid in
$C(K)$.  
\end{abstract}
\maketitle

\section{Introduction}
\label{sec:introduction}
The classical approximation theorem due to P. P.~Korovkin \cite{Kor53}
states that the asymptotic behaviour of any sequence of positive linear
maps on $C[0,1]$ is completely determined by its behaviour on the three
functions $\boldsymbol{1}$, $x$ and $x^2$. The underlying reason for this
phenomenon was revealed by Yu.~A.~\v{S}a\v{s}kin \cite{Sas67}, who
characterized all such test sets in terms of the Choquet boundary. To
state his result, let $X$ be a compact metric space and let
$G\subseteq C(X)$ be a set whose linear span $S$ contains the constant
function $\boldsymbol{1}$. The \emph{Choquet boundary} $\partial_S X$ of
$S$ is the set of all points $t\in X$ such that the Dirac measure
$\delta_t$ is the only Borel probability measure $\mu$ on $X$ satisfying
$f(t)=\int_X f\,d\mu$ for all $f\in S$. \v{S}a\v{s}kin's theorem asserts
that $\partial_S X=X$ if and only if for any sequence of positive linear
maps $\varPhi_k\colon C(X)\to C(X)$ ($k\in \natu$),
   \begin{align*}
\lim_{k\to\infty}\|\varPhi_k(g)-g\|=0 \;\; \forall g\in G \implies
\lim_{k\to\infty}\|\varPhi_k(f)-f\|=0 \;\; \forall f\in C(X).
   \end{align*}
In other words, the asymptotic behaviour of the sequence
$\{\varPhi_k\}_{k=1}^\infty$ on the $C^*$-algebra $C(X)$ is uniquely
determined by the set $G$ precisely when every point of $X$ belongs to
the Choquet boundary of $S$.

\v{S}a\v{s}kin's theorem thus shows that the Choquet boundary lies at the
heart of classical approximation theory. Inspired by this observation,
and by the central role that approximation plays in the modern theory of
operator algebras, W.~Arveson \cite{Arv11} introduced the notion of
\emph{hyperrigidity}, a ``noncommutative'' counterpart of the Korovkin
property that turns out to encompass a wide range of significant
operator-algebraic phenomena.

For the reader's convenience, we recall Arveson's definition (see
\cite[Definition~1.1]{Arv11}). Here and below, UCP stands for unital
completely positive.
\begin{df} \label{dyfnh}
A finite or countably infinite set $G$ of generators of a unital
$C^*$-algebra $\ascr$ is said to be \textit{hyperrigid} if for any
faithful representation $\pi\colon \ascr\to \ogr(\hh)$ on a Hilbert
space $\hh$ and for any sequence $\varPhi_k\colon \ogr(\hh)\to \ogr(\hh)$
($k\in \natu$) of UCP~maps,
   \begin{align*}
\lim_{k\to\infty}\|\varPhi_k(\pi(g))-\pi(g)\|=0 \;
\forall g\in G \implies
\lim_{k\to\infty}\|\varPhi_k(\pi(a))-\pi(a)\|=0 \;
\forall a\in \ascr.
   \end{align*}
\end{df}

In the spirit of \v{S}a\v{s}kin's theorem, Arveson
\cite[Conjecture~4.3]{Arv11} conjectured that a generating set is
hyperrigid whenever every irreducible representation of the generated
$C^*$-algebra is a boundary representation. The converse implication is
always valid. Bilich and Scherer \cite{BilichScherer} recently constructed
counterexamples to Arveson's hyperrigidity conjecture in commutative
$C^*$-algebras, showing that the boundary condition does not imply the
unique extension property for arbitrary representations.

Let $I\subset\rbb$ be a compact interval and let $f\in C(I)$ be
real-valued and strictly convex. The function system
$S_f=\operatorname{span}\{\boldsymbol{1},x,f\}$ has Choquet boundary $I$
(see \cite[Propositions~9.1 and~9.2]{Arv11}).
Arveson \cite[Remark~9.5]{Arv11} asked whether $\{I,T,f(T)\}$ is
hyperrigid in the unital $C^*$-algebra generated by $T$ whenever $T$ is
selfadjoint with spectrum in $I$. Brown \cite{Brown} answered this
question affirmatively as a consequence of the following weak-to-strong
convergence theorem.
Below, $\sigma(T)$ denotes the spectrum of an operator $T$.
   \begin{theorem}[{\cite[Theorem~2.1]{Brown}}]
\label{brosowo} Let $f$ be a strictly convex
continuous function on an interval $J\subseteq \rbb$,
let $\hh$ be a Hilbert space, and let $\{T_i\}_{i\in
\varOmega}$ be a net of selfadjoint operators in
$\ogr(\hh)$ such that $\sigma(T_i)\subseteq J$ for all
$i\in \varOmega$. If $\{T_i\}_{i\in \varOmega}$
converges weakly to a selfadjoint operator $T\in
\ogr(\hh)$ with $\sigma(T)\subseteq J$, and if also
$\{f(T_i)\}_{i\in \varOmega}$ converges weakly to
$f(T)$, then $\{\varphi(T_i)\}_{i\in \varOmega}$
converges strongly to $\varphi(T)$ for every bounded
continuous function $\varphi$ on $J$. In particular if
the net $\{T_i\}_{i\in \varOmega}$ is bounded, then
$\{T_i\}_{i\in \varOmega}$ converges strongly to $T$.
   \end{theorem}

The present paper extends this rigidity phenomenon to strictly convex
functions of two variables, under a regularity assumption. 

\begin{theorem}
\label{main1}
Let $K\subset\rbb^2$ be a nonempty compact convex set, and let
$f\colon K\to\rbb$ be a strictly convex function admitting a
$C^2$ extension to an open neighbourhood of $K$.
Then the set $G=\{x,y,f\}$ generates $C(K)$ and is hyperrigid in
$C(K)$.  
\end{theorem}

\section{Preliminaries}
\label{sec:preliminaries}

Throughout the paper, $\rbb$, $\comp$ and $\natu$ denote the
sets of real numbers, complex numbers and positive integers,
respectively. For a topological space $X$, we write $\borel(X)$
for the $\sigma$-algebra of its Borel subsets. The symbol
$\boldsymbol{1}_D$ stands for the characteristic function of
a set $D$, and $\lin Z$ denotes the complex linear span of
a subset $Z$ of a complex vector space.

Let $\hh$ be a complex Hilbert space. We denote by $\ogr(\hh)$
the $C^*$-algebra of all bounded linear operators on $\hh$,
with identity operator $I=I_{\hh}$. An operator
$T\in\ogr(\hh)$ is positive if
$\langle Th,h\rangle\Ge0$ for every $h\in\hh$.

\section{Proof of Theorem~\ref{main1}}
\label{sec:brown-proof}

Choose an open neighbourhood $U\subset\rbb^2$ of $K$ and
a function $\widetilde f\in C^2(U,\rbb)$ such that
$\widetilde f|_K=f$. We denote by
$D^2\widetilde f(u)
=(\partial_i\partial_j\widetilde f(u))_{i,j=1}^2$
the Hessian matrix of $\widetilde f$ at $u\in U$.

For each $a\in K$, define $\psi_a\colon K\to\rbb$ by
$\psi_a(t)=f(t)-f(a)
-\langle\nabla\widetilde f(a),t-a\rangle$.
Strict convexity implies that $\psi_a(t)\Ge0$,
with equality if and only if $t=a$.
Taylor's formula with integral remainder yields
\begin{equation}\label{eq:quadratic-gap}
\begin{aligned}
0\Le\psi_a(t)
&=\int_0^1(1-s)
\bigl\langle
D^2\widetilde f(a+s(t-a))(t-a),t-a
\bigr\rangle\,\D s\\
&\Le L\abs{t-a}^2,\qquad a,t\in K,
\end{aligned}
\end{equation}
where $L>0$ satisfies
$L\Ge\frac12\sup_{u\in K}\norm{D^2\widetilde f(u)}$.

Strict convexity of $f$ implies that $\psi_a(t)\Ge0$ for
all $t\in K$, with equality if and only if $t=a$.

\begin{lemma}\label{lem:localization}
Let $K\subset\rbb^2$ be nonempty and compact, and let
$E\colon\borel(K)\to\ogr(\hh)$ be a spectral measure. Set
\begin{align}
G_z=\int_K \abs{t-z}^2\,\D E(t),\qquad z\in\rbb^2.
\end{align}
If a positive operator $T$ satisfies
\begin{align}
T\Le cG_a,\qquad a\in K,
\end{align}
for some constant $c>0$, then $T=0$.
\end{lemma}

\begin{proof}

Choose a square $Q_0\subset\rbb^2$ of side length $L>0$ such that
$K\subset Q_0$. We shall repeatedly subdivide each square into four
congruent subsquares by bisecting its sides. At every stage, boundary
points are assigned to exactly one of the adjacent subsquares. We refer
to the resulting Borel sets as \emph{dyadic cells}. Thus, each dyadic
cell is the disjoint union of its four children and is contained in a
closed square of the corresponding side length.

Fix $h\in\hh$ and set $\eta=T^{1/2}h$. We may assume that
$\eta\ne0$. Define $\mu\colon\borel(\rbb^2)\to[0,\infty)$ by
\begin{align*}
\mu(B)
&=\big\langle E(K\cap B)\eta,\eta\big\rangle
 =\big\|E(K\cap B)\eta\big\|^2,
 \qquad B\in\borel(\rbb^2).
\end{align*}
Then $\mu$ is a finite positive Borel measure on $\rbb^2$, supported
on $K$, and
\begin{align*}
\mu(\rbb^2)=\|\eta\|^2.
\end{align*}

Fix $N\geq1$. Among the four children of $Q_0$, choose a dyadic cell
$Q_1$ satisfying $\mu(Q_1)\geq\mu(Q_0)/4$, and continue inductively
within the cell so obtained. This yields a nested sequence of dyadic
cells
\begin{align*}
Q_0\supset Q_1\supset\cdots\supset Q_N
\end{align*}
such that, for $0\Le j\Le N$,
\begin{align*}
\mu(Q_j)&\Ge 4^{-j}\|\eta\|^2,\\
\operatorname{diam}(Q_j)^2&\Le 2L^2 4^{-j}.
\end{align*}
Since $\mu(Q_N)>0$ and $\mu$ is supported on $K$, we have
$K\cap Q_N\ne\varnothing$. Fix
\begin{align*}
a\in K\cap Q_N.
\end{align*}

For $0\Le j\Le N$, set
\begin{align*}
P_j=E(K\cap Q_j),\qquad
D_N=\sum_{j=0}^N4^jP_j,
\qquad
s_N=\langle D_N\eta,\eta\rangle.
\end{align*}
By the definition of $\mu$ and the choice of the cells $Q_j$,
\begin{align}\label{eq:dyadic-local-mass}
s_N
 =\sum_{j=0}^N4^j\mu(Q_j)
 \Ge (N+1)\|\eta\|^2>0.
\end{align}
Define $d_N\colon K\to[0,\infty)$ by
\begin{align*}
d_N(t)=\sum_{j=0}^N4^j\mathbf{1}_{Q_j}(t),
\qquad t\in K.
\end{align*}
Then
\begin{align*}
D_N=\int_K d_N(t)\,\D E(t).
\end{align*}
Fix $t\in K$, and let $m\in\{0,\ldots,N\}$ be the largest index
such that $t\in Q_m$. Since the cells are nested and
$a\in Q_N\subset Q_m$, we have
\begin{align*}
d_N(t)
&=\sum_{j=0}^m4^j\Le\frac43\,4^m,\\
|t-a|^2
&\Le\operatorname{diam}(Q_m)^2\Le2L^2 4^{-m}.
\end{align*}
Consequently,
\begin{align}\label{eq:dyadic-pointwise}
|t-a|^2d_N(t)\Le\frac83L^2,\qquad t\in K.
\end{align}
Multiplying \eqref{eq:dyadic-pointwise} by $d_N(t)$ and applying
the spectral calculus, we obtain
\begin{align}\label{eq:dyadic-operator}
D_NG_aD_N
=\int_K|t-a|^2d_N(t)^2\,\D E(t)
\Le\frac83L^2D_N.
\end{align}

Since $\eta=T^{1/2}h$, the Cauchy--Schwarz inequality,
the assumption $T\Le cG_a$, and \eqref{eq:dyadic-operator} give
\begin{align*}
s_N^2
&=\bigl|\langle D_N\eta,T^{1/2}h\rangle\bigr|^2\\
&=\bigl|\langle T^{1/2}D_N\eta,h\rangle\bigr|^2\\
&\Le\|h\|^2\langle TD_N\eta,D_N\eta\rangle\\
&\Le c\|h\|^2\langle D_NG_aD_N\eta,\eta\rangle\\
&\Le\frac{8cL^2}{3}\|h\|^2s_N.
\end{align*}
As $s_N>0$, division by $s_N$ and
\eqref{eq:dyadic-local-mass} yield
\begin{align}\label{eq:dyadic-local-bound}
(N+1)\|T^{1/2}h\|^2
\Le\frac{8cL^2}{3}\|h\|^2.
\end{align}
Although the selected cells and the point $a$ may depend on $h$
and $N$, the constant on the right-hand side is independent
of both. The estimate also holds when $T^{1/2}h=0$.
Since $h\in\hh$ is arbitrary, it follows that
$\|T\|\Le8cL^2/[3(N+1)]$ for every $N\geq1$.
Letting $N\to\infty$ gives $T=0$.
\end{proof}

\begin{proof}[Proof of Theorem~\ref{main1}]
The coordinate functions separate the points of $K$, so
$\{\boldsymbol{1},x_1,x_2,f\}$ generates $C(K)$. By the unique
extension characterization of hyperrigidity, it suffices to show
that any unital representation
$\pi\colon C(K)\to\ogr(\hh)$ and any UCP map
$\Phi\colon C(K)\to\ogr(\hh)$ satisfying
\[
\Phi(s)=\pi(s),\qquad s\in\{\boldsymbol{1},x_1,x_2,f\},
\]
coincide on $C(K)$.
Let $E$ be the spectral measure of $\pi$, and let $F$ be the
normalized semispectral measure representing $\Phi$.

For $a\in K$, set
\[
\psi_a(t)
=f(t)-f(a)-\langle\nabla f(a),t-a\rangle,
\qquad t\in K.
\]
Strict convexity and Taylor's formula give
\begin{align}
\label{eq:uep-supporting-gap}
0\Le\psi_a(t)\Le L|t-a|^2,
\qquad
\psi_a(t)=0\ \Longleftrightarrow\ t=a,
\qquad a,t\in K,
\end{align}
where $L>0$ is independent of $a$ and $t$.
Moreover, since
$\psi_a\in\operatorname{span}\{\boldsymbol{1},x_1,x_2,f\}$,
\begin{align}
\label{eq:uep-gap-moments}
\int_K\psi_a\,\D F
=\Phi(\psi_a)
=\pi(\psi_a)
=\int_K\psi_a\,\D E.
\end{align}

Let $A,B\subset K$ be nonempty disjoint compact sets.
Joint continuity of $(a,t)\mapsto\psi_a(t)$ and
\eqref{eq:uep-supporting-gap} imply that
\[
m:=\min\{\psi_a(t):a\in A,\ t\in B\}>0.
\]
For every $a\in A$, the inequality
$m\boldsymbol{1}_B\Le\psi_a$ holds on $K$. Integrating with
respect to $F$, using \eqref{eq:uep-gap-moments}, and compressing
to $E(A)\hh$, we obtain
\begin{align}
\label{eq:uep-local-domination}
0\Le E(A)F(B)E(A)
&\Le \frac{1}{m}\int_A\psi_a(t)\,\D E(t) \notag\\
&\Le \frac{L}{m}\int_A|t-a|^2\,\D E(t),
\qquad a\in A.
\end{align}
If $E(A)\ne0$, the restriction of $E$ to $A$, regarded as a
spectral measure on $E(A)\hh$, is normalized. Thus
Lemma~\ref{lem:localization} applied to
\eqref{eq:uep-local-domination} yields
$E(A)F(B)E(A)=0$. This equality is immediate if $E(A)=0$.
Positivity of $F(B)$ now gives
\begin{align}
\label{eq:uep-separated-sets}
E(A)F(B)=F(B)E(A)=0.
\end{align}

Fix a nonempty compact set $B\subset K$ and put
\[
A_n=\{t\in K:\dist(t,B)\Ge 1/n\},
\qquad n\in\natu.
\]
The sets $A_n$ are compact, disjoint from $B$, and increase to
$K\setminus B$. By \eqref{eq:uep-separated-sets} and the strong
countable additivity of $E$,
\[
E(K\setminus B)F(B)=F(B)E(K\setminus B)=0.
\]
Since $0\Le F(B)\Le I$, it follows that
\begin{align}
\label{eq:uep-compact-domination}
F(B)=E(B)F(B)E(B)\Le E(B).
\end{align}
The same inequality holds trivially for $B=\varnothing$.

Let $D\in\borel(K)$ and $h\in\hh$. Regularity of the finite
positive Borel measure
$D\mapsto\langle F(D)h,h\rangle$ and
\eqref{eq:uep-compact-domination} imply that
\[
\langle F(D)h,h\rangle
=
\sup_{\substack{B\subset D\\ B\text{ compact}}}
\langle F(B)h,h\rangle
\Le \langle E(D)h,h\rangle.
\]
Consequently, $F(D)\Le E(D)$ for every $D\in\borel(K)$.
Applying this inequality to $K\setminus D$ and using
$F(K)=E(K)=I$, we obtain the reverse inequality. Hence $F=E$,
and integration gives $\Phi=\pi$.
\end{proof}

\end{document}